\documentclass[12pt]{amsart}
\usepackage[]{amsmath, amsthm, amsfonts, verbatim, amssymb, tikz}

\usepackage{amsmath}
\usepackage{amsthm}
\usepackage{amssymb}
\usepackage{hyperref}
\usepackage{rotating}
\usepackage{amssymb}
\usepackage{epstopdf}
\usepackage{pifont}
\usepackage{amsmath}
\usepackage{enumerate}
\usepackage{graphicx}
\DeclareMathAlphabet{\mathpzc}{OT1}{pzc}{m}{it}
\usepackage{tikz}
\usepackage{tikz-cd}
\usepackage{cases}
\usetikzlibrary{arrows}
\usepackage{adjustbox}

\newtheorem{theorem}{Theorem}[section]
\newtheorem{lemma}[theorem]{Lemma}
\newtheorem{proposition}[theorem]{Proposition}
\newtheorem{corollary}[theorem]{Corollary}

\theoremstyle{remark}

\newcommand{\bB}{\mathbb{B}}
\newcommand{\bC}{\mathbb{C}}

\newcommand{\bQ}{\mathbb{Q}}
\newcommand{\bP}{\mathbb{P}}
\newcommand{\bR}{\mathbb{R}}
\newcommand{\bZ}{\mathbb{Z}}

\newcommand{\cC}{\mathcal{C}}

\newcommand{\cG}{\mathcal{G}}

\newcommand{\cL}{\mathcal{L}}
\newcommand{\cO}{\mathcal{O}}

\newcommand{\cT}{\mathcal{T}}

\newcommand{\Aut}{{\rm Aut}}

\newcommand{\Gal}{{\rm Gal}}

\newcommand{\NS}{{\rm NS}}

\newcommand{\Pic}{{\rm Pic}}

\newcommand{\oGamma}{\overline \Gamma}

\def\<{\langle}
\def\>{\rangle}

\def\ni{\noindent}
\def\ra{\rightarrow}

\def\tt{\widetilde{t}}

\title[Surfaces of maximal canonical degree]
 {Examples of Surfaces with Canonical Map of Maximal Degree} 

\author{Ching-Jui Lai}
\address{Department of Mathematics, National Cheng Kung University, Tainan 70101\\ Taiwan (R.O.C.)}
\email{cjlai72@mail.ncku.edu.tw}

\author{Sai-Kee Yeung}
\address{Mathematics Department, Purdue University, West Lafayette, IN  47907\\USA }
\email{yeung@math.purdue.edu}

\begin{document}

\begin{abstract} It was shown by A. Beauville that if the canonical map $\varphi_{|K_M|}$ of a complex smooth projective surface $M$ is generically finite, then $\deg(\varphi_{|K_M|})\leq 36$. The first example of a surface with canonical degree 36 was found by the second author. In this article, we show that for any surface which is a degree four Galois \'etale cover of a fake projective plane $X$ with the largest possible automorphism group $\Aut(X)=C_7:C_3$ (the unique non-abelian group of order 21), the base locus of the canonical map is finite, and we verify that 35 of these surfaces have maximal canonical degree 36. We also classify all smooth degree four Galois \'etale covers of fake projective planes, which give possible candidates for surfaces of canonical degree $36$. Finally, we also confirm in this paper the optimal upper bound of the canonical degree of smooth threefolds of general type with sufficiently large geometric genus, related to earlier work of C. Hacon and J.-X. Cai.\end{abstract}

\maketitle

\section{Introduction}

Let $M$ be a smooth complex projective minimal surface of general type with $p_g(M)\neq0$. Assume that the canonical map, $$\varphi=\varphi_{|K_M|}:M\dashrightarrow W:=\overline{\varphi(M)}\subseteq\bP^{p_g(M)-1}$$
is generically finite onto its image. We are interested in the \emph{canonical degree} of $M$, the degree of $\varphi$. If $\varphi$ is not generically finite, we simply say that $M$ has canonical degree zero. The following proposition was proved in \cite{B}, cf. \cite{Y1}. We include the proof here for completeness. 
\begin{proposition}\label{Bea} Let $M$ be a minimal surface of general type whose canonical map $\varphi=\varphi_{|K_M|}$ is generically finite. Then $\deg \varphi\leq 36$. Moreover, $\deg\varphi=36$ if and only if $M$ is a smooth ball quotient $\bB_\bC^2/\Sigma$ with $p_g(M)=3$, $q(M)=0$, and $|K_M|$ is base point free.
\end{proposition}
\begin{proof} Let $P$ be the mobile part of $|K_M|$.  Let  $S\rightarrow M$ be a resolution of $P$ and let $P_S$ be the induced base point free linear system defining $S\rightarrow W=\overline{\varphi(M)}$. Then 
\begin{align*} \deg\varphi\cdot(p_g-2)\leq\deg\varphi\cdot\deg W
                                                            = P_S^2\leq P^2\leq K_M^2
                                                            \leq9\chi(\cO_M)\leq9(1+p_g).
\end{align*}
The first inequality is the degree bound for a non-degenerate surface in $\bP^n$ given in \cite{B}, while 
the fourth inequality is the Bogomolov-Miyaoka-Yau inequality. Hence as $p_g\geq3$, we have 
\begin{align*} \deg\varphi\leq 9(\frac{1+p_g}{p_g-2})\leq36.
\end{align*}
Moreover, $\deg\varphi=36$ only when $p_g(M)=3$, $q(M)=0$, and $P_S^2=P^2=K_M^2$. This is only possible when $|K_M|$ is base point free.  In such a case,
 $K_M^2=36=9\chi(\cO_M)$ and hence $M$ is a smooth ball quotient $\bB_\bC^2/\Sigma$ by results of Aubin and Yau, cf. \cite{B} or \cite{BHPV}.
\end{proof}

\noindent{\bf Notation.} Throughout this paper, we do not distinguish line bundles with divisors. The linear equivalence and numerical equivalence of divisors are written respectively as $D_1\sim D_2$ and $D_1\equiv D_2$. The cyclic group of order $n$ is denoted by $C_n$. The group $C_7:C_3$ is the unique non-abelian group of order 21. The projective space of dimension $n$ over $\bC$ is denoted by $\bP^n$. A finite field of order $n$ is denoted by $F_n.$

\vskip 0.2 cm

From Proposition \ref{Bea}, it is an interesting problem to know the geometric realization of possible canonical degrees and many surfaces with canonical degree at most $16$ have been constructed, see \cite{P} or \cite{DG} for more references. However, the first example of a surface with maximal canonical degree 36 was constructed only recently by \cite{Y1} as a suitably chosen $C_2\times C_2$-Galois cover of a special fake projective plane $X$. The fake projective plane $X$ in \cite{Y1} has ${\rm Aut}(X)=C_7:C_3$, and by \cite{LY} it satisfies $h^0(X,2L_X)=0$ for every ample generator $L_X$ of ${\rm NS}(X)$. The choice of the lattice for the ball quotient $M$ is explicitly described in \cite{Y1} via the classifying data of \cite{PY} and \cite{CS}. 

Here are the main goals of this paper. The first goal is to construct more examples of surfaces with maximal canonical degree. This is given as Theorem \ref{main} below.  Then we examine the corresponding question in complex dimension 3, given as Corollary \ref{3fold} below. A second goal is to identify all potential examples of surfaces of canonical degree 36 constructed as a degree four Galois \'etale cover of a fake projective plane. We prove that for these Galois covers the canonical maps have at worst discrete base locus whenever the underlying fake projective plane has the largest possible automorphism group $C_7:C_3$. This is given as Theorem \ref{iso} and Proposition \ref{all}. For the presentation of this paper, we start with Theorem \ref{iso} hoping that it would give the reader a more comprehensible overall picture. 

We remark that our proof of Theorem \ref{main} is essentially independent of Theorem \ref{iso} and Proposition \ref{all}.  A reader who is interested only in new surfaces of canonical degree $36$ may briefly go over statements in earlier sections and proceed directly to Section \ref{new} of the paper.

Recall that a fake projective plane is a ball quotient $X=\bB_\bC^2/\Pi$ for some lattice $\Pi\subseteq{\rm PU}(2,1)$, where $\Pi$ is constructed as a subgroup of a maximal arithmetic lattice $\oGamma$. An unramified cover $M$ of $X$ is given by 
$\bB^2_\bC/\Sigma$ for a normal subgroup $\Sigma\lhd\Pi$ of finite index. For the sequence of Galois covers 
$$M:=B_{\bC}^2/\Sigma\stackrel{p}\rightarrow X=B_{\bC}^2/\Pi\stackrel{q}\rightarrow B_{\bC}^2/\oGamma$$
corresponding to the normal subgroups $\Sigma\lhd\Pi\lhd\oGamma$, one has the covering group ${\rm Gal}(M/X)=\Pi/\Sigma$ and $\Aut(X)=\oGamma/\Pi$. We focus on the case when $|{\rm Gal}(M/X)|=4$ and $\Aut(X)=C_7:C_3$. Our first theorem identifies potential examples of surfaces of canonical degree 36.

\begin{theorem}\label{iso} Let $M\rightarrow X$ be a degree four Galois \'etale cover over a fake projective plane $X$ with $\Aut(X)=C_7:C_3$. Then $q(M)=0$ and the base locus of the linear system $|K_M|$ is discrete.
\end{theorem}

A degree four Galois \'etale cover $M\rightarrow X$ over a fake projective plane $X$ is determined by a quotient of $H_1(X,\bZ)$ of order four, to be explained in details in Lemma \ref{cover} of Section \ref{pre}. The degree of this cover is dictated by the possible existence of a surface of maximal canonical degree, i.e., $K_M^2/K_X^2=4$. There are many degree four covers of fake projective planes. For future reference, we classify all such surfaces. In the table below, only lattices of fake projective planes giving rise to Galois \'etale covers of degree four are listed, which is the case if there is a normal subgroup of index four in the lattice $\Pi$ corresponding to a given fake projective plane $X=\bB^2_\bC/\Pi$. This list of the fake projective planes follows the conventions in \cite{PY} and \cite{CS}. 
In the following table, we have 
\begin{enumerate}
	\item column 1: $k$ is a totally real number field, $\ell$ is a totally imaginary extension of $k$, and $\mathcal T$ represents a finite number of places relevant to the classification. These are notations used to classify fake projective planes defined in \cite{PY};
	\item column 2 the corresponding naming of classes of maximal arithmetic lattices containing fake projective planes in \cite{CS} corresponding to $\overline{\Gamma}$ in the notation of \cite{PY}, where $a$ and $p$ are data from the first column;
	\item column 3: the naming of the individual fake projective planes in each class used in \cite{CS};
	\item column 4: $\mbox{Aut}(X)$ is the automorphism group of a fake projective plane $X$;
	\item column 5: the first homology class of a fake projective plane $X$;
	\item column 6: $N_0$ is the number of degree $4$ coverings of $X$, which is the number of subgroups of index four of the lattice $\Pi$;
	\item column 7: $N_1$ denotes the number of normal coverings among the degree $4$ coverings above.
\end{enumerate}
All the examples in the last column satisfy $H_1(M,\mathbb{Q})=0$, which implies $q(M)=0$ by Poincar\' e Duality.

\vskip 0.2 cm

\begin{adjustbox}{center, rotate=0, nofloat, caption=[Table 1]}
$\begin{array}{|c|c|c|c|c|c|c|}
\hline
(k,\ell,\cT)&\mbox{class}&X&{\rm Aut}(X)&H_1(X,\bZ)&N_0&N_1\\ \hline\hline
(\bQ,\bQ(\sqrt{-1}),\{5\})&(a=1,p=5,\emptyset)&(a=1,p=5,\emptyset, D_3)&C_3&C_2\times C_4\times C_{31}&4&3\\ \cline{2-7}
&(a=1,p=5,\{2\})&(a=1,p=5, \{2\},D_3)&C_3&C_4\times C_{31}&4&1\\ \cline{1-7}
(\bQ,\bQ(\sqrt{-1}),\{2,5\})&(a=1,p=5,\{2I\})&(a=1,p=5,\{2I\})&\{1\}&C_2\times C_3\times C_4^2&47&19\\ \hline
(\bQ,\bQ(\sqrt{-2}),\{3\})&(a=2,p=3,\emptyset)&(a=2,p=3,\emptyset, D_3)&C_3&C_2^2\times C_{13}&4&1\\ \cline{2-7}
&(a=2,p=3,\{2\})&(a=2,p=3, \{2\},D_3))&C_3&C_2^2\times C_{13}&4&1\\ \cline{1-7}
(\bQ,\bQ(\sqrt{-2}),\{2,3\})&(a=2,p=3,\{2I\})&(a=2,p=3,\{2I\})&\{1\}&C_2^4\times C_3&83&35\\ \hline
(\bQ,\bQ(\sqrt{-7}),\{2\})&(a=7,p=2,\emptyset) &(a=7,p=2,\emptyset, D_3 2_7)&C_7:C_3&C_2^4&91&35\\ \cline{3-7}
&&(a=7,p=2,\emptyset,7_{21})&\{1\}&C_2^2\times C_3\times C_7&3&1\\ \cline{2-7}
&(a=7,p=2,\{7\})&(a=7,p=2,\{7\},D_3 2_7)&C_7:C_3&C_2^3&7&7\\  \cline{3-7}
&&(a=7,p=2,\{7\},D_3 7'_7)&C_3&C_2^2\times C_7&2&1\\ \cline{3-7}
&&(a=7,p=2,\{7\},7_{21})&\{1\}&C_2^3\times C_3&19&7\\ \hline
(\bQ,\bQ(\sqrt{-7}),\{2,3\})&(a=7,p=2,\{3\})&(a=7,p=2,\{3\},D_3)&C_3&C_2\times C_4\times C_7&4&3\\ \cline{3-7}
&&(a=7,p=2,\{3\},3_3)&\{1\}&C_2^2\times C_3\times C_4&19&11\\ \cline{2-7}
&(a=7,p=2,\{3,7\})&(a=7,p=2,\{3,7\},D_3)&C_3&C_4\times C_7&2&1\\ \cline{3-7}
&&(a=7,p=2,\{3,7\},3_3)&\{1\}&C_2\times C_3\times C_4&7&3\\ \hline
(\bQ,\bQ(\sqrt{-7}),\{2,5\})&(a=7,p=2,\{5\})&(a=7,p=2,\{5\})&\{1\}&C_2^2\times C_9&3&1\\ \cline{1-7}
(\bQ,\bQ(\sqrt{-15}),\{2\})&(a=15,p=2,\emptyset)&(a=15,p=2,\emptyset,D_3)&C_3&C_2^2\times C_7&2&1\\  \cline{3-7}
&&(a=15,p=2,\emptyset,3_3)&\{1\}&C_2^3\times C_9&11&7\\ \cline{2-7}
&(a=15,p=2,\{3\})&(a=15,p=2,\{3\},3_3)&C_3&C_2^3\times C_3&19&7\\ \cline{2-7}
&(a=15,p=2,\{5\})&(a=15,p=2,\{5\},3_3)&\{1\}&C_2^2\times C_9&3&1\\ \cline{2-7}
&(a=15,p=2,\{3,5\})&(a=15,p=2,\{3,5\},3_3)&C_3&C_2^2\times C_3&1&1\\ \hline
(\cC_{18},\{v_3\})&(\cC_{18},p=3,\emptyset)&(\cC_{18},p=3,\emptyset,d_3 D_3)&C_3\times C_3&C_2^2\times C_{13}&1&1\\  \hline
(\cC_{20},\{v_2\})&(\cC_{20},\{v_2\},\emptyset)&(\cC_{20},\{v_2\},\emptyset,D_3 2_7)&C_7:C_3&C_2^6&651&651\\  \cline{2-7}
&(\cC_{20},\{v_2\},\{3+\})&(\cC_{20},\{v_2\},\{3+\},D_3)&C_3&C_4\times C_7&2&1\\ \cline{3-7}
&&(\cC_{20},\{v_2\},\{3+\},\{3+\}_3)&\{1\}&C_2\times C_3\times C_4&7&3\\  \cline{2-7}
&(\cC_{20},\{v_2\},\{3-\})&(\cC_{20},\{v_2\},\{3-\},D_3)&C_3&C_4\times C_7&2&1\\ \cline{3-7}
&&(\cC_{20},\{v_2\},\{3-\},\{3-\}_3)&\{1\}&C_2\times C_3\times C_4&7&3 \\ \hline
\end{array}$
\end{adjustbox}
\vskip 0.2cm
\begin{center}
{\sc Table 1}
\end{center}

\vspace{-\baselineskip} 
\begin{proposition}\label{all} There are altogether $835$ lattices which give rise to $1670$ non-biholomorphic smooth minimal surfaces as degree four Galois \'etale covers of fake projective planes with $q(M)=0$. 
\end{proposition}

From Table 1, there are 35 degree four Galois \'etale covers of the fake projective plane $(a=7, p=2,\emptyset, D_32_7)$, which all have Galois group $C_2\times C_2$. Generalizing the result of \cite{Y1}, we show that these \'etale covers all have canonical degree 36.

 \begin{theorem}\label{main} The $35$ degree four  Galois \'etale covers of the fake projective plane $(a=7, p=2,\emptyset, D_32_7)$, all with Galois group $C_2\times C_2$, are minimal surfaces of general type with canonical degree $36$.
\end{theorem}

Our result has the implication on the optimal canonical degree for smooth threefolds of general type with large geometric genus. We refer the readers to Section \ref{sec3fold} for more details.

\begin{corollary}\label{3fold} There exist many examples of smooth minimal threefolds of general type $Y$ with the degree of the canonical map $\deg(\Phi_{|K_Y|})=72$.  In fact, there exist such threefolds with $p_g(Y)=3g$ and $K_Y^3=72(g-1)$ for each $g\geqslant 2$.
 \end{corollary}

The surface studied in \cite{Y1} has Picard number one, which is a deep result in automorphic forms from \cite{Ro}, \cite{BR}, and is used in \cite{Y1} to simplify the geometric arguments. For a general degree four \'etale cover of a fake projective plane, it is not clear whether the Picard number equals to one. Comparing to the result in \cite{Y1}, one technical improvement in the present article is to show that any surface as in Theorem \ref{iso} possesses a generically finite canonical map. Continuing from this, mobility of the canonical system is proved but in a different argument from \cite{Y1}. In fact, we can show that any degree four \'etale cover of a fake projective plane with $\Aut(X)=C_7:C_3$ has generically finite canonical map and at worst discrete base locus. To get rid of the finite number of base points, we need more detailed information about the canonical sections as given in \cite{Y1}, see in particular the corrigendum there. By analyzing carefully the method used in \cite{Y1}, we come up with new examples of surfaces with maximal canonical degree by considering new degree four Galois \'etale covers of the same fake projective plane $X$ used in \cite{Y1}.  These new \'etale covers correspond to various $C_2\times C_2$ quotient groups of $H_1(X,\bZ)=C_2^4$. In such cases, we are able to write down relevant global sections explicitly with the help of Magma and finish the prove of base point freeness. This last step is where we have to restrict further the type of lattice $\Sigma$ associated to $M$.
 
To find which \'etale cover works for our scheme, as a first step we list all normal subgroups of index four in a lattice associated to a fake projective plane. All fake projective planes supporting such a subgroup are listed in the third column of Table 1 above. Now for each of the listed surfaces, we exhaust all possible normal subgroups of index four. The procedure of finding such a surface as well as verification of necessary conditions stated in Theorem \ref{iso} and Proposition \ref{all} is similar to that in \cite{Y1}. In \cite{Y1}, the choice of the $C_2\times C_2$ Galois \'etale cover is very specific and has to come from killing the $2$-torsion invariant line bundles under a Sylow $3$-subgroup of the automorphism group $C_7:C_3$. In this paper, we obtain more examples by overcoming this technical hurdle, namely, we consider all possible $C_2\times C_2$ Galois \'etale covers of the fake projective plane in \cite{Y1}.

The explicit computation is accomplished by using Magma. The proof of Theorem \ref{main} generalizes the argument of \cite{Y1}. 

Here is the organization of this paper. We first prepare some preliminary results related to our construction in Section \ref{pre}. The proofs of Theorem \ref{iso} and \ref{main} are given in Section \ref{seciso} and \ref{new} respectively. 
Finally we study the corresponding problem in dimension three in Section \ref{sec3fold}.

\section{Preliminary discussions and idea of proofs}\label{pre}
Let $X=\bB_\bC/\Pi$ be a fake projective plane with $\pi_1(X)=\Pi$. It is known from definition that the first Betti number of $X$ is trivial.  According to \cite{PY}, there is always a nontrivial 
torsion element in $H_1(X,\bZ)$. The torsion group $H_1(X,\bZ)$ is available from \cite{CS}.

\begin{lemma}\label{cover} A fake projective plane $X$ possesses a degree four Galois \'etale cover if and only if there is a quotient group of order four of $H_1(X,\bZ)$. 
\end{lemma}
\begin{proof} We know that $H_1(X,\bZ)$ is a direct sum of finite cyclic abelian groups as the first Betti number of $X$ is trivial. If $Q$ is a quotient group of order four of $H_1(X,\bZ)$, then there is a homomorphism
$$\rho:\Pi\rightarrow\Pi/[\Pi,\Pi]=H_1(X,\bZ)\rightarrow Q.$$
The kernel of $\rho$ gives rise to a normal subgroup $\Sigma$ of index four in $\Pi$, with $Q$ as the deck transformation group of the covering map $M=\bB^2_\bC/\Sigma\rightarrow X=\bB^2_\bC/\Pi$.

On the other hand, if there is a normal subgroup $\Sigma$ of index four in $\Pi$, it leads to a homomorphism $\sigma:\Pi\rightarrow \Pi/\Sigma$. As a group of order four is always abelian, $\sigma$ factors through a homomorphism $\Pi/[\Pi,\Pi]\rightarrow\Pi/\Sigma$.  We conclude that $\Pi/\Sigma$ lives as a quotient group of order four of $\Pi/[\Pi,\Pi]=H_1(X,\bZ)$.
\end{proof}

We consider an \'etale cover $\pi:M\rightarrow X$ corresponding to a subgroup $\pi_1(M)\leq\Pi$ of index four. In particular, the finite group $\cG=\Pi/\pi_1(M)$ is either $C_2\times C_2$ or $C_4$. 

\begin{lemma}\label{gg} Let $M$ be a smooth projective surface and assume that there is an \'etale  cover $\pi:M\rightarrow X$ of degree four over a fake projective plane $X$. Suppose that $q(M)=0$, then $p_g(M)=3$. 
\end{lemma}
\begin{proof} Since $\pi:M\rightarrow X$ is \'etale and $p_g(X)=q(X)=0$, $\chi(\cO_M)=4\chi(\cO_X)=4$. It follows that $p_g(M)=3$ if $q(M)=0$. 
\end{proof}

Suppose now a surface $M$ is constructed as in Lemma \ref{gg}. We study 
the canonical map $\varphi=\varphi_{|K_M|}:M\dashrightarrow\bP^2$. We will assume that $\pi:M\rightarrow X$ is a \emph{Galois cover}, i.e., $\Sigma:=\pi_1(M)\leq\Pi$ is normal. 
Note that then $|K_M|$ is invariant under the Galois group $\cG:={\rm Gal}(M/X)=\Pi/\pi_1(M)$.

Let us relate the canonical sections from Lemma \ref{gg} to divisors on $X$.
It is known from the Universal Coefficient Theorem that torsions in $H_1(X.\bZ)$ give rise to a torsion line bundle on $X$, cf. Lemma 4 of \cite{LY}. Denote by $\cL_\chi$ the invertible sheaf on $X$ corresponding to a torsion line bundle on $X$ given by a character $\chi$.  In this case, the trivial character $\cO_X$ is denoted by $\cL_1$.
The push forward of the structure sheaf of $M$ splits into eigen-sheaves
$$
\pi_*\cO_M=\bigoplus_{\chi:\cG\rightarrow \bC^*}\cL_\chi,
$$
Denote by $\omega_M$ the dualizing sheaf of a surface $M$.  Then
$$\pi_*\omega_M=\bigoplus_{\chi:\cG\rightarrow \bC^*}\omega_X\otimes \cL_\chi.$$
It follows from the degeneration of the Leray spectral sequence that 
\begin{equation}
H^i(M,\omega_M)=\bigoplus_{\chi:\cG\rightarrow \bC^*}H^i(X,\omega_X\otimes \cL_\chi)
\end{equation}
for all $i$.  Hence vanishing of $q(M)$ implies that $H^1(X,\omega_X\otimes \cL_\chi)=0$ for
all $\chi:\cG\rightarrow \bC^*$.  By Serre Duality, $h^2(X,\omega_X\otimes\cL_\chi)=h^0(X,\cL_\chi^{-1})$, which is either
$0$ or $1$ depending on whether $\chi$ is trivial of not.  From Riemann-Roch formula and the fact that $X$ is a fake projective plane,
it follows that $h^0(X,\omega_X\otimes\cL_\chi)=1$ for each $\chi\neq 1$, which corresponds to three linearly independent sections in Lemma \ref{gg}.
Denote by $D_1, D_2, D_3$ the corresponding curves on $X$.  It follows that $H^0(M,K_M)$ is generated by $\pi^*(D_i)$, $i=1,2,3$, noting that
$\pi^* \cL_\chi\cong\cO_M$.

\begin{lemma}\label{bpf23} Assume that $q(M)=0$ and let $D_1, D_2, D_3$ be divisors obtained as above.  Assume that $D_1\cap D_2\cap D_3=\emptyset$.  Then
$H^0(M,K_M)$ is base point free and the canonical degree of $M$ is $36$.
\end{lemma}

\begin{proof} Let $x$ be a point in the base point set of $|K_M|$.  Since $|K_M|$ is invariant under the Galois group $\cG$, $\pi(x)\in D_1\cap D_2\cap D_3$, which
is empty.  It follows from Proposition \ref{Bea} that the canonical degree of $M$ is $36$.
\end{proof}

The last lemma would be utilized in Section \ref{new} to give a proof of Theorem \ref{main}.  The presentation here is a simplification of the original one,
thanks to the suggestion of the referee.

\section{General constraints on base point set}\label{seciso}

The goal of this section is to give a proof of Theorem \ref{iso}, which gives constraints on the base point set of $|K_M|$ without knowledge on an explicit description of fake projective plane $X$. Here as $\rho(X)=1$, we always denote by $L_X$ an ample generator of $\Pic(X)$. Also recall that for a fake projective plane $X$, we have $p_g(X)=q(X)=0$ and $L_X^2=1$ by definition. We begin with the following simple observations. 

\begin{lemma}\label{gen} Let $X$ be a fake projective plane and let $L_X$ be an ample generator of $\Pic(X)$. Then $h^0(X,L)\leq1$ for any line bundle $L\equiv L_X$ and $h^0(X,L')\leq 2$ for any line 
	bundle $L'\equiv 2L_X$.
	\end{lemma}
\begin{proof} If $L''$ is a line bundle with $L''\equiv 4L_X$, then by Riemann-Roch formula $h^0(X,L'')=3$. But if $L\equiv L_X$ and $H^0(X,L)$ has two linearly independent sections $x$ and $y$, then 
	$\{x^4,x^3y,x^2y^2,xy^3,y^4\}$ are five linearly independent sections of $H^0(X,L^{\otimes 4})$, which is absurd. The second statement is proved similarly. 
	\end{proof}

\begin{lemma}\label{Sch} If $C$ is an irreducible and reduced curve on a fake projective plane $X$ with $C\equiv L_X$, then $C$ is smooth of genus 3. 
\end{lemma}
\ni{\bf Proof.} Given an irreducible and reduced curve $C$, we denote by $C^\nu$ the normalization of $C$ and $\nu:C^\nu\rightarrow C$ the normalization morphism. The $\cO_C$ sheaf
$\delta:=\nu_*\cO_{C^\nu}/\cO_C$ is the cokernel of the natural map $\cO_C\rightarrow\nu_*\cO_{C^\nu}$ and satisfies 
$$g(C^\nu)=p_a(C^\nu)=p_a(C)-h^0(C,\delta).$$

We first remark that $g(C^\nu)\geq2$ as $X$ is hyperbolic. The Ahlfors-Schwarz Lemma applied to the composition map induced by the normalization $\nu':C^\nu\xrightarrow{\nu}C\hookrightarrow X$ (cf. \cite{CCL}) for the manifolds equipped with Poincar\'e metrics implies that the K\"ahler forms satisfy $\nu'^*\omega_X\leq\omega_{C^\nu}$, with equality if and only if it is a holomorphic isometry leading to totally geodesic $C$. Since there is no totally geodesic curve on a fake projective plane from the proof of \cite[Lemma 6]{LY}, the inequality is strict. Hence for $C\equiv kL_X$ with $k\geq1$, integrating over $C^\nu$, we get
$$2k= \frac{2}{3}(K_X\cdot C)<\deg(K_{C^\nu})=2g(C^\nu)-2=k(k+3)-2h^0(C,\delta),$$
where we used the fact that the Ricci curvature is $\frac{3}{2}$ of the holomorphic sectional 
curvature for the Poincar\'e metric on $X$ and the adjunction $p_a(C)=\frac{1}{2}C\cdot(K_X+C)$. Hence $k = 1$ implies that $h^0(C,\delta) = 0$ and $C$ is smooth with $g(C)=3$. 
\qed

\begin{lemma}\label{inv} Let $X$ be a fake projective plane with a nontrivial automorphism group and let $C$ be an effective divisor such that $C\equiv L_X$. For any nontrivial subgroup $H\leq\Aut(X)$ with $H\cong C_3$ or $C_7$, $h^*C\neq C$ for any $h\in H-\{e\}$.
	\end{lemma}
\begin{proof} Clearly $C$ must be reduced and irreducible as $\rho(X)=1$. From Lemma \ref{Sch}, $C$ is smooth of genus three. Suppose now $h^*C=C$ for all $h\in H$. From \cite[Lemma 6]{LY}, $H$ must act non-trivially on $C$. Note that $H$ can only be $C_3$ or $C_7$ from the list of \cite{CS}. 

If $H\cong C_7$, then there exists an $H$-fixed point on $C$, as by the Hurwitz formula there is no \'etale cover of degree 7 from a smooth genus three curve. 
By \cite[Lemma 7]{LY}, for $x=\dim_\bC H^1(C,\cO_C)^{\rm inv}$ we have the equation,
$$n=2-2\cdot3+\frac{2\cdot7}{7-1}(3-x)\ \Rightarrow\ 3n+7x=9. $$
The only solution is $(n,x)=(3,0)$ and $C/C_7\subseteq X/C_7$ is a smooth rational curve. But then there is a non-constant lifted map from $\bP^1$ to the universal cover $\bB^2_\bC$ of $X/C_7$, this contradicts to Liouville's theorem. 

If $H\cong C_3$, then there exists an $H$-fixed point on $C$, as by the Hurwitz formula there is no \'etale cover of degree 3 from a smooth genus three curve. By the same argument as above, we see that $(n,x)=(5,0)$ or $(2,1)$. In either cases, there is a non-constant lifted map from $\bP^1$ or $\bC$ to $\bB^2_\bC$, which again contradicts Liouville's theorem.
	\end{proof}	

\begin{lemma}\label{propgenfin} Let $X$ be a fake projective plane with $\Aut(X)=C_7:C_3$. Suppose that there is a Galois \'etale cover $\pi:M\rightarrow X$ of degree four and $q(M)=0$, then the canonical map $\varphi:M\dashrightarrow\bP^2$ is generically finite. 
\end{lemma}
\begin{proof} From Lemma \ref{gg}, we know that $p_g(M)=3$ and hence the canonical map maps $M$ to $\bP^2$. Write $|K_M|=P+F$, where $P$ is the mobile part and $F$ is the fixed divisor.  By construction, we have $\varphi=\varphi_{|K_M|}=\varphi_{P}:M\dashrightarrow\bP^2$.  We will  abuse the notation: $P$ will be the mobile linear system or a general member in it. 

Assume that $\overline{\varphi(M)}=C\subseteq\bP^2$ is a curve. We will derive a contradiction.

First of all, we claim that $P$ is not base point free, or equivalently $P^2\neq0$. Assume now $P^2=0$. We consider $\cG=\Gal(M/X)$. Since $g^*K_M=K_M$ for any $g\in\cG$, we have that $g^*F=F$ for each $g\in\cG$. Indeed, $g^*P$ is a mobile sub-linear system of $|K_M|$ and hence $g^*F\geq F$ as Weil divisors. Hence as $\pi$ is Galois, $F=\pi^*F_X$ for an effective divisor $F_X$ on $X$. Moreover, if $\NS(X)=\<L_X\>$ for an ample divisor $L_X$, then $K_X\equiv3L_X$, $F_X\equiv lL_X$ for some $0\leq l\leq 3$, and $P\equiv\pi^*(3-l)L_X$. Now, $P^2=0$ implies that $l=3$ and hence $P\equiv0$. This is a contradiction as a non-zero effective divisor cannot be numerically trivial. 

Since $\varphi:M\dashrightarrow C\subseteq\bP^2$ is not a morphism, we take a composition of finitely many smooth blow-ups $\rho:\widehat{M}\rightarrow M$ to resolve $P$ and let $\psi:\widehat{M}\rightarrow C\subseteq\bP^2$ be the induced morphism. We have the following diagram after taking the Stein factorization of $\psi:S\rightarrow C$:
\begin{center}
\begin{tikzpicture}
implies/.style={double double equal sign distance, -implies},
\node (m) at (0,2) {$\widehat{M}$};
\node (M) at (0,0) {$M$};
\node (C) at (2,0) {$C$};
\node (in) at (2.5,0) {$\subseteq$};
\node (P) at (3,0) {$\bP^2$};
\node (CC) at (2,2) {$\tilde{C}$};

\path[->] (m) edge node[left]{$\rho$}(M);
\path[->] (m) edge node[above]{$\beta$} (CC);
\path[dashed,->] (M) edge node[below]{$\varphi$} (C);
\path[->] (CC) edge node[right]{$\alpha$} (C);
\path[->] (m) edge node[above]{$\psi$} (C);
\end{tikzpicture}
\end{center}

If $\rho^*P=\widehat{P}+\widehat{F}$, where $\widehat{P}=\psi^*|\cO_C(1)|$ is base point free, $\widehat{F}\geq0$ is the fixed divisor, and $\psi=\psi_{\widehat{P}}$, then $\widehat{F}$ is a non-trivial effective $\rho$-exceptional divisor with $\beta(\widehat{F})=\tilde{C}$. In particular, $\tilde{C}\cong\bP^1$ as all the irreducible components of $\widehat{F}$ are rational. Since $\alpha:\tilde{C}\rightarrow C$ is defined by $\alpha^*|\cO_C(1)|\subseteq|\cO_{\bP^1}(d)|$ for some $d\geq1$ and hence an element in $\widehat{P}$ is given by $\beta^*H$ for some $H\in|\cO_{\bP^1}(d)|$, we have $\widehat{P}\supseteq\beta^*|\cO_{\bP^1}(d)|$. In particular, we get  
$$\widehat{P}=\psi^*|\cO_C(1)|=\beta^*\alpha^*|\cO_C(1)|=\beta^*|\cO_{\bP^1}(d)|.$$ 
As $\dim \widehat{P}=p_g(M)=3$, we get $d=2$ and $C\subseteq\bP^2$ being irreducible and non-degenerate is a smooth conic in $\bP^2$. 

Let $\widehat{M}_{c}$ be a general fibre of $\widehat{M}\rightarrow\tilde{C}$ and $D:=\rho_*(\widehat{M}_{c})\equiv P/2$ be the corresponding prime divisor on $M$. Recall that $\pi:M\rightarrow X$ is Galois, $K_M=\pi^*K_X\equiv\pi^*(3L_X)$ and $P\equiv\pi^*(lL_X)$ for some $1\leq l\leq 3$ as $P^2\neq0$, where $\NS(X)=\langle L_X\rangle$ and $L_X^2=1$. It follows from the genus formula, 
$$(K_M+D)\cdot D=2g_a(D)-2\in2\bZ$$ 
that $l=2$ is the only possibility. Hence $P\equiv\pi^*(2L_X)$, $F=\pi^*F_X\equiv\pi^*L_X$, and $D\equiv\pi^*L_X$. Note that if $h^0(X,2L_X)=0$ for any ample generator $L_X$ on $X$, then we arrive the required contradiction as $2F_X\neq0$. This is exactly the argument in \cite{Y1}, where the vanishing holds for $X$ a very special fake projective plane as discussed in the introduction. Below we provide a more elementary argument. 

It is easy to see that $\cG$ acts on $C\cong\bP^1$ holomorphically and induces an action on $\tilde{C}$. We claim that there is always a fixed point on $\tilde{C}=\bP^1$. If $\cG$ acts trivially, then every point is a fixed point.\footnote{In fact, this case is absurd. If $\cG$ acts trivially on $C$, then $\cG$ also acts trivially on $\tilde{C}\cong\bP^1$. Any fibre of $\beta:\widehat{M}\ra\tilde{C}$ as a section of $H^0(\bP^1,\cO_{\bP^1}(1))$ is $\cG$-fixed and descends to a $\cG$-invariant section $D\equiv\pi^*L_X$ on $M$, which then descends to a section $D_X\equiv L_X$ on $X$. For any two such sections $D$ and $D'$ on $M$, $D\sim D'$ implies that $D_X\equiv D'_X\equiv L_X$ where $\pi^*D_X=D$ and $\pi^*(D_X')=D'$.  Since $X$ has only finitely many nontrivial torsion but $H^0(\bP^1,\cO_{\bP^1}(1))$ is infinite, we can find a line bundle $L=L_X+T_X$ for some torsion line bundle $T_X$ on $X$ with $\dim|L|\geq1$. This contradicts Lemma \ref{gen}.} Otherwise, $\cG$ has two fixed points on $\tilde{C}$ from the Lefschetz fixed point formula. In particular, the fiber $\widehat{M}_c$ over a fixed point $c$ is $\cG$-invariant and descends to an effective divisor $G^X\equiv L_X$ on $X$.\footnote{Up to here everything works for all fake projective planes with a nontrivial automorphism group.} 	 

Suppose now that $\Aut(X)=C_7:C_3$. Note that in this case a non-trivial torsion elements is always a 2-torsion. In particular for any $\sigma\in\Aut(X)$, 
$\sigma^*G^X\sim G^X+T_\sigma$ for some 2-torsion $T_\sigma$ and 
$$\sigma^*(2G^X)=2\sigma^*(G^X)\sim 2G^X+2T_\sigma=2G^X.$$ 
On the other hand, for any non-trivial element $\sigma\in\Aut(X)$, $G^X\neq\sigma^*G^X$ by Lemma \ref{inv}. The curves $G^X$ and $\sigma^*G^X$ intersect at a unique point $Q_\sigma$ as $G^X\cdot(\sigma^*G^X)=L_X^2=1$. We claim that there are three linearly independent sections of the form $2\sigma^*G^X$ in $|2G^X|$, which then contradicts to Lemma \ref{gen}.

We fix one non-trivial $\sigma$ and consider $Q:=Q_\sigma$. Note that then 
$2G^X$ intersects with $\sigma^*(2G^X)$ only at $Q$ with multiplicity four. By the result of \cite{PY}, the isotropic group at $Q$ cannot be the whole $\Aut(X)$. Hence there exists a nontrivial element $\tau\in\Aut(X)$, $\tau\neq\sigma$, such that $\tau^*Q\neq Q$. 
In particular, $\tau^*(2G^X)$ only intersects with $\tau^*\sigma^*(2G^X)$ at $\tau^*Q$ with multiplicity four. Since elements in the pencil $\<\mu\cdot2G^X+\lambda\cdot2\sigma^*G^X\>$ must pass through $Q$ with multiplicity four, one of $\tau^*(2G^X)$ and $\tau^*\sigma^*(2G^X)$ is not in $\<2G^X,2\sigma^*G^X\>$ or otherwise $\tau^*Q=Q$. Hence $h^0(X,2G^X)>2$ and we have a contradiction to Lemma \ref{gen}. 

Hence we conclude that $\dim\overline{\varphi(M)}\neq1$. Since $\varphi(M)\subseteq\bP^2$ has to be positive dimensional, we conclude that $\varphi:M\dashrightarrow\bP^2$ must be dominant and hence generically finite.
\end{proof}

\begin{lemma}\label{codim1} Let $M\rightarrow X$ be a Galois \'etale cover of degree four of a fake projective plane $X$ with $\Aut(X)=C_7:C_3$. If $q(M)=0$, then the canonical linear system $|K_M|$ is mobile, i.e., there is no codimension one base locus.
\end{lemma}
\begin{proof} We follow the same notation as in the proof of Lemma \ref{propgenfin}: $\NS(X)=\<L_X\>$ for an ample divisor $L_X$, $K_X\equiv3L_X$, $F_X\equiv lL_X$ for some $0\leq l\leq 3$, and $P\equiv\pi^*(3-l)L_X$. We claim that $l=0$. 

Since $\dim P=p_g(M)-1=2>0$, $P$ contains a nontrivial effective divisor and hence $l\neq 3$. 

If $l=1$, then we consider the action of $\Aut(X)=C_7:C_3$ on $F_X=L_X+T$, where $T$ is a 2-torsion. Then the same argument as in the proof of Lemma \ref{propgenfin} produces a line bundle $\cL\equiv 2L_X$ with $h^0(X,\cL)>2$, but this violates Lemma \ref{gen}. 

If $l=2$, then we consider the same argument as above on $P_X\equiv L_X$. 

 Here is an alternate argument.  In the above setting, if $H^0(X,2L_X)=0$ for $L_X$ any ample generator of $\Pic(X)$, then $|K_M|=P$ is mobile. Indeed, the assumption also implies that $H^0(X,L_X)=0$ for any ample generator of $\Pic(X).$ Hence for $F=\pi^*F_X$ with $F_X\equiv lL_X$, $l=0$ is the only possibility and $F=0$. The hypothesis holds for any fake projective plane with an automorphism group of order 21 by a result of \cite{LY}. 
\end{proof}

{\it Proof of Theorem \ref{iso}} First of all, from Magma, all Galois coverings of a fake projective plane of index $4$ can be listed, as is done in the proof of Proposition \ref{all} below.
Furthermore, Magma tells us that abelianization of the lattices associated to such coverings are all trivial.  Hence $q(M)=0$ for our examples. Theorem \ref{iso} now follows from Lemma \ref{codim1}.

\begin{proof}[of Proposition \ref{all}] We simply apply the procedure of construction as in \cite{Y1} to each of the fake projective plane listed in column 3 of Table 1. We first need to enumerate all possible surfaces as degree four Galos \'etale cover associated to fake projective planes as listed. It turns out  that the number of index four subgroups of the lattice $\Pi$ to a fake projective plane
in the table is recorded in the column $N_1$ in Table 1. This could be seen by considering subgroups of order $4$ in $H_1(X,\bZ)$ as in Lemma \ref{cover}, or by listing index four subgroups of $\Pi$ from Magma.

Now we claim that all the different sub-lattices of index $4$ of $\Pi$ in Table 1 give rise to non-isometric complex hyperbolic forms in terms of the Killing metrics on the locally symmetric spaces. For this purpose, we assume that $\Lambda_1$ and $\Lambda_2$ are two groups obtained from the above procedure and $B_{\bC}^2/\Lambda_1$ is isometric to $B_{\bC}^2/\Lambda_2$. From construction, $\Lambda_1$ and $\Lambda_2$ are normal subgroups of index 4 in two lattices $\Pi_1$ and $\Pi_2$ corresponding to the fundamental groups
of fake projective planes. Let $\oGamma_1$ and $\oGamma_2$ be the corresponding maximal arithmetic groups in the respective classes. As $B_{\bC}^2/\Lambda_1$ and  $B_{\bC}^2/\Lambda_2$ are isometric, $\Lambda_1$ is conjugate to $\Lambda_2$ as discrete subgroups of the same algebraic group $G$ with $G\otimes \bR\cong PU(2,1)$.  Hence the two corresponding maximal lattices satisfy $\oGamma_1\cong \oGamma_2$, and similarly $\Pi_1\cong\Pi_2$. It follows that they have to come from the same row in the Table 1 and hence correspond to the same subgroup of index $4$ in the same lattice associated to some fake projective plane.  Hence there are altogether $835$ non-isometric complex two ball quotients obtained in this way, by summing over the column of $N_1$ in Table 1.

Now for each locally symmetric space $M=B_{\bC}^2/\Lambda$ obtained as above, it gives rise to a pair of complex structures $J_1$ and $ J_2$, which are conjugate to each other.
These two complex structures give rise to two non-biholomorphic complex surfaces $S_1=(M,J_1)$ and $S_2=(M,J_2)$.  In fact, if they are biholomorphic, the corresponding
four-fold quotient $S_1/[\Pi,\Lambda]$ and $S_2/[\Pi,\Lambda]$ are biholomorphic and are fake projective space. This contradicts the results in \cite{KK}, see also the Addendum of \cite{PY}, that conjugate complex structures on a fake projective space give rise to two different complex structures.  
 
 In general, let  $(M_1, J_1)$ and $(M_2,J_2)$ be two complex ball quotients obtained from taking degree 4 \'etale covers of some possibly different fake projective planes.  If $(M_1, J_1)$ and $(M_2,J_2)$ are biholomorphic, they are isometric with respect to the corresponding Bergman (Killing) metrics.  Hence from the earlier argument, $M_1$ is isometric to $M_2$ and we may regard $M_1=M_2$.  Now the argument of the last paragraph implies that $J_1=J_2$.  In conclusion, we conclude that the $1670$ complex surfaces  obtained from the pair of conjugate complex structures on the $835$ underlying locally symmetric structures give rise to distinct complex surfaces. This concludes the proof of Proposition \ref{all}.
\end{proof}

\section{New examples of surfaces with maximal canonical degree}\label{new}
 Our goal in this section is to prove Theorem \ref{main}. The surface studied in \cite{Y1} and here is constructed from the fake projective plane $X$ given in \cite[Section 5.9]{PY} in the class of $(a=7, p=2)$ and is denoted by $(a=7, p=2,\emptyset, D_32_7)$ in the notation of \cite{CS}.

\begin{proof}[of Theorem \ref{main}] We consider $\pi:M\rightarrow X$ a Galois $C_2\times C_2$--\'etale  cover of the fake projective plane $X$ in the class $(a=7, p=2,\emptyset, D_32_7)$. From Magma computation, the irregularity $q(M)=0$, cf. Proposition \ref{all}. Hence by Lemma \ref{bpf23}, it suffices for us to prove that the canonical map of $M$ is base point free. From the discussion in Section \ref{pre}, there are non-trivial 2-torsions $\tau_i\in\Pic^0(X)$ for $i=1,2,3$ corresponding to characters of $\cG=\Gal(M/X)=C_2\times C_2$ such that $H^0(X,K_X+\tau_i)=\<t_i\>$ and $H^0(M,K_M)=\<\pi^*t_i|\ i=1,2,3\>$.

For the convenience of the reader, we recall the key steps of the argument in \cite{Y1}. For simplicity, we denote by $G$ the automorphism group $\Aut(X)=C_7:C_3$. The automorphism group of $X$ has a presentation 
$G=\langle a,b|a^7=b^3=1, bab^{-1}=a^2\rangle.$ The group $G$ contains a normal Sylow 7-subgroup $G_7=\<a\>$, and seven conjugate Sylow $3$-subgroups, one of which is $G_3:=\<b\>$.  We know from the Riemann-Roch formula that $h^0(X,2K_X)=10$. In terms of the explicit basis of $H^0(X,2K_X)$ given by \cite{BK}, the action of $G$ is presented by 
\begin{eqnarray}
&&a(u_0:u_1:u_2:u_3:u_4:u_5:u_6:u_7:u_8:u_9)\nonumber \\
&=&(u_0:\zeta_7^6u_1:\zeta_7^5u_2:\zeta_7^3u_3:\zeta_7u_4:\zeta_7^2u_5:\zeta_7^4u_6:\zeta_7u_7:\zeta_7^2u_8:\zeta_7^4u_9)\\
&&b(u_0:u_1:u_2:u_3:u_4:u_5:u_6:u_7:u_8:u_9) \nonumber \\
&=&(u_0:u_2:u_3:u_1:u_5:u_6:u_4:u_8:u_9:u_7)
\end{eqnarray}

From the Corrigendum of \cite{Y1}, under the action of $G_7$, $S:=\cup_{\Sigma\in C_2^4-\{1\}} H^0(X,K_X+\Sigma)$ consists of 3 orbits, where we recall that a $p$-torsion element $\Sigma\in H_1(X,\bZ)=C_2^4$ correspond to a $p$-torsion element $\Sigma\in\Pic^0(X)$ by the universal coefficient theorem (see \cite[Lemma 4]{LY}). 
    	\begin{enumerate}
    		\item $\<\tt_0\>=H^0(X,K_X+\Sigma_0)$, where $\Sigma_0$ is $G$-invariant corresponding to an element in $H_1(X/G,\bZ)^\times$ and $\tt_0^2=u_0$. 
    		\item Two disjoint $G_7$ orbits $\<a\>\tt_1$ and $\<a\>\tt_2$, where $\tt_i$'s are $G_3$-invariant corresponding to elements in $H_1(X/G_3,\bZ)^\times-\{\Sigma_0\}$. 		
 	\end{enumerate}
Let $v_0=u_0, v_1=u_1+u_2+u_3, v_2=u_4+u_5+u_6$, and $v_3=u_7+u_8+u_9$. From \cite{Y1}, one finds that \begin{equation}\label{4.3}\begin{cases} \tt_0^2=v_0,\\ \tt_1^2=v_0+\frac12(1+\sqrt{-7})v_1,\\ \tt_2^2=v_0+(-5+\sqrt{-7})v_1+4(1-\sqrt{-7})v_2-4(v_3)\end{cases} 
\end{equation}
with the help of elementary command \verb'IsDomain' in Magma. It is proved that $\cap_{i=0}^2Z_{t_i^2}=\emptyset$, which was verified in the Corrigendum of \cite{Y1} by checking that $\cap_{i=0}^2Z_{t_i^2}=\emptyset$ on $X$ modulo $p=23$ from the command \verb'HilbertPolynomial' in Magma. We remark that the same example was also studied later in \cite{Ri}, where the author independently verified with more sophisticated techniques in Magma that the sections obtained from the above procedure do give rise to sections in $H^0(M,K_M)$.

Now under the action of $G_7$, the explicit sections $\tt_0$ and $a^j\tt_i$, $i=1,2$ and $0\leq j\leq 6$, precisely give the effective sections of $S:=\cup_{\Sigma\in C_2^4-\{1\}} H^0(X,K+\Sigma)$. We will prove that $\cap_{i=0}^2Z_{t^2_i}=\emptyset$ by consider possible choices of $\{t_1,t_2,t_3\}\subseteq S=\<\tt_0\>\cup\<a\>\tt_1\cup\<a\>\tt_2$ and check by Magma whether these sections have common intersection.

Conjugating by an element in $G_7$, we may assume that $t_1$ belongs to $\{\tt_0,\tt_1,\tt_2\}$.  Suppose $t_1=\tt_0$, where $\tt_0$ is invariant as a set under $G$, then conjugate by an element in $G_7$, we may assume that $t_2=\tt_1$. But by construction $\tau_3=\tau_1\cdot\tau_2$ is determined by $\tau_1=\sigma_0$ and $\tau_2$, which gives $\tt_2\in H^0(X,K_X+\tau_0\cdot\tau_1)=H^0(X,K_X+\tau_2)$. In particular, this case was already checked in \cite{Y1} as 
$\cap_{i=0}^2Z_{t^2_i}= Z_{v_0}\cap Z_{\tt_1^2}\cap Z_{\tt^2_2}=\emptyset$
and we are done. 

Consider now the case that none of $t_i$'s is $\tt_0$.  In this scenario, $t_i$ belongs to the orbits of $\tt_1$ or $\tt_2$.  Again we use the fact that effective divisors $D_i$'s have common intersections if and only if $2D_i$'s have common intersections. Hence it suffices for us to prove the following claim.
\begin{lemma}\label{check} Let $i, j\in \{1,\dots,6\}$. Then
\begin{enumerate}[$(a)$]
	\item  $Z_{\tt_1}\cap Z_{a^i\tt_1}\cap Z_{a^j\tt_2}=\emptyset$ for $1\leq i,j\leq6;$
	\item $Z_{\tt_2}\cap Z_{a^i\tt_1}\cap Z_{a^j\tt_2}=\emptyset$ for $1\leq i,j\leq6;$
	\item  $Z_{\tt_1}\cap Z_{a^i\tt_1}\cap Z_{a^j\tt_1}=\emptyset$ for $1\leq i<j\leq6;$
	\item $Z_{\tt_2}\cap Z_{a^i\tt_2}\cap Z_{a^j\tt_2}=\emptyset$ for $1\leq i<j\leq6.$
\end{enumerate}
\end{lemma}

\begin{proof} In terms of the basis chosen with action of $G_7$ given in equation (4.1) and the explicit sections listed in (4.3), statement $(a)$ in Lemma \ref{check} holds if there is no the common intersection for the following sections,
\begin{eqnarray*}
&\{&u_0+\frac12(1+\sqrt{-7})(u_1+u_2+u_3),\\
&&u_0+\frac12(1+\sqrt{-7})(\zeta_7^{-i}u_1+\zeta_7^{-2i}u_2+\zeta_7^{-4i}u_3),\\
&&u_0+(-5+\sqrt{-7})(\zeta_7^{-j}u_1+\zeta_7^{-2j}u_2+\zeta_7^{-4j}u_3)+4(1-\sqrt{-7})(\zeta_7^{j}u_4
+\zeta_7^{2j}u_5+\zeta_7^{4j}u_6)\\&&-4(\zeta_7^{j}u_7+\zeta_7^{2j}u_8+\zeta_7^{4j}u_9)\}.
\end{eqnarray*}
 Instead of using the command \verb'HilbertPolynomial' over the cyclotomic field $\bQ(\zeta_7)$ on $X$, we specialize it to the finite field $F_{29}$, where $16$ is a primitive $7$-th root of unity and $14$ serves as $\sqrt{-7}$.  In this way, computing over the finite field $F_{29}$, we verify from Magma that the above three polynomials do not have common intersection on $X$ for all $i, j\in \{1,\dots,6\}$ in $F_{29}$. This implies that the original
equations do not have common zero over the algebraic number field $\bQ(\zeta_7)$. Similar arguments applies to $(b)$, $(c)$, and $(d)$ in the Lemma \ref{check}.
\end{proof}

We remark that Lemma \ref{check} actually is stronger than what is sufficient for our purpose. For example, consider the case of $(a)$. It is enough to check $Z_{\tt_1}\cap Z_{a^i\tt_1}\cap Z_{a^j\tt_2}=\emptyset$ for one pair of $(i,j)$ corresponding to the elements $\cG-\{1\}=\{\tau_1,\tau_2,\tau_3\}.$ However, since
we are checking by Magma, the extra computation does not make any essential difference in computer time. Similar argument applies to the cases $(b)$, $(c)$, $(d)$ as well.

Theorem \ref{main} follows immediately from Lemma \ref{check}.
\end{proof}

\section{Remark on maximal canonical degree of threefolds}\label{sec3fold}
Theorem \ref{main} has an implication on the canonical degree bound of threefolds. The purpose of this section is to explain literatures in this direction and relations to Theorem \ref{main}. From this point on, let $Y$ be a Gorenstein minimal complex projective threefold of general type with locally factorial terminal singularities. Suppose that the linear system $|K_Y|$ defines a generically finite map $\Phi=\Phi_{|K_Y|}:Y\dashrightarrow\bP^{p_g(Y)-1}.$ M. Chen asked in \cite{Ch} if there is an upper bound of $\deg(\Phi)$. A positive answer was provided in \cite{Hac} with $\deg(\Phi)\leq576.$ Later on, it was improved in \cite{DG2} that   $\deg(\Phi)\leq360$ (with equality if and only if $p_g(Y)=4, q(Y)=2, \chi(\omega_Y)=5, K_Y^3=360$, and $|K_Y|$ is base point free.) In \cite{C}, it is shown that $\deg(\Phi)\leq 72$ if the geometric genus satisfies $p_g(Y)>10541$. 

As a corollary of Theorem \ref{main} and the above discussion, we conclude that the canonical degree 72  can be achieved as stated in Corollary \ref{3fold}.
 
 \begin{proof}[of Corollary \ref{3fold}] Equipped with Theorem \ref{main}, the corollary follows essentially from an observation of \cite[Section 3]{C}.  

Take $C$ a smooth hyperelliptic curve of genus $g\geq2$, then the canonical map $\varphi_{|K_C|}:C\rightarrow\bP^{g-1}$ is the composition of the double cover $C\rightarrow\bP^1$ with the $(g-1)$-Veronese embedding $\bP^1\hookrightarrow\bP^{g-1}$. In particular, $\deg(\varphi_{|K_C|})=2$, cf. \cite{Har}. Take $M$ a surface satisfying the optimal degree bound $\deg(\varphi_{|K_M|})=36$ as in Theorem \ref{main}, then $\varphi=\varphi_{|K_M|}:M\rightarrow\bP^2$ is a generically finite morphism of $\deg(\varphi)=K_M^2=36$.

Now let $Y=X\times C$, then $Y$ is a smooth projective threefold of general type with $p_g(Y)=3g$ and $\Phi=\Phi_{|K_Y|}:Y\rightarrow\bP^{3g-1}$ a morphism. From our construction, it follows that $\Phi$ is generically finite and   
$$\deg{\Phi}\cdot\deg W=K_Y^3=3K_X^2\cdot K_C=3\cdot36\cdot(2g-2),$$ 
where $W=\Phi(Y)$ is the image of the composition maps  $Y\hookrightarrow\bP^2\times\bP^{g-1}\hookrightarrow\bP^{3g-1}$ defined by $|K_Y|$ and ${\cO_{\bP^2\times\bP^{g-1}}(1,1)}$. Hence $\deg W=3(g-1)$ and $\deg(\Phi)=72.$
\end{proof}

\ni{\bf Acknowledgements.}\label{ackref}
It is a pleasure for the second author to thank Donald Cartwright for his help on Magma commands. The authors would like to express their appreciation and thankfulness to the referee for very helpful comments and suggestions on the paper. This work is partially done during the first author's visit at Research Institute of Mathematical Sciences in Kyoto, National Center of Theoretical Sciences and National Taiwan University in Taiwan, and the second author's visit of the Institute of Mathematics of the University of Hong Kong.  The authors thank the warm hospitality of the institutes.

\end{document}